\documentclass{amsart}

\usepackage{mathtools,amsmath,amsthm,amsfonts,latexsym,amssymb,mathrsfs,setspace,verbatim,cite,mathtools,tikz-cd,appendix,comment,enumitem}
\usepackage{fullpage}

\mathtoolsset{showonlyrefs=true}

\newtheorem{thm}{Theorem}

\newtheorem{lem}[thm]{Lemma}
\newtheorem{prop}[thm]{Proposition}

\numberwithin{thm}{section}
\numberwithin{equation}{section}

\newtheorem{question}[thm]{Question}

\newcommand{\rat}{\mathbb Q}
\newcommand{\real}{\mathbb R}

\newcommand{\alg}{\overline\rat}
\newcommand{\algt}{\alg^{\times}}

\newcommand{\intg}{\mathbb Z}

\newcommand{\tors}{\mathrm{tors}}

\newcommand{\supp}{\mathrm{supp}}
\newcommand{\spann}{\mathrm{span}}

\newcommand{\xx}{{\bf x}}
\newcommand{\zz}{{\bf z}}

\newcommand{\lcc}{LC_c(Y_S)}
\newcommand{\lcco}{LC_c^0(Y_S)}

\begin{document}

\title[Consistent Maps]{Consistent Maps and Their Associated Dual Representation Theorems}

\author[C.L. Samuels]{Charles L. Samuels}
\address{Christopher Newport University, Department of Mathematics, 1 Avenue of the Arts, Newport News, VA 23606, USA}
\email{charles.samuels@cnu.edu}

\begin{abstract}
	A 2009 article of Allcock and Vaaler examined the vector space $\mathcal G := \overline{\mathbb Q}^\times/\overline{\mathbb Q}^\times_{\mathrm{tors}}$ over $\mathbb Q$, describing its completion with respect to the Weil 
	height as a certain $L^1$ space.  By involving an object called a consistent map, the author began efforts to establish Riesz-type representation theorems for the duals of spaces related to $\mathcal G$.  Specifically, we provided such results
	for the algebraic and continuous duals of $\overline{\mathbb Q}^\times/{\overline{\mathbb Z}}^\times$.  In the present article, we use consistent maps to provide representation theorems for the duals of locally 
	constant function spaces on the places of $\overline{\mathbb Q}$ that arise in the work of Allcock and Vaaler. 
	We further apply our new results to recover,  as a corollary, a main theorem of our previous work.
\end{abstract}

\maketitle

\section{Introduction}

\subsection{Background} Let $\alg$ be a fixed algebraic closure of $\rat$ and let $\algt_\tors$ denote the group of roots of unity in $\algt$.  Following the notation of \cite{AllcockVaaler}, we write $\mathcal G = \algt/\algt_\tors$
and note that $\mathcal G$ is a vector space over $\rat$ with addition and scalar multiplication given by
\begin{equation} \label{GOperations}
	(\alpha,\beta) \mapsto \alpha\beta\quad\mbox{and}\quad (r,\alpha) \mapsto \alpha^r.
\end{equation}
For each number field $K$, we write $M_K$ to denote the set of all places of $K$.   If $L/K$ is a finite extension and $w\in M_L$, then $w$ divides a unique place $v$ of $K$, and in this case, 
we shall write $K_w$ to denote the completion of $K$ with respect to $v$.  Additionally, we let $p_v$ be the unique place of $\rat$ such that $v$ divides $p_v$, so in our notation we have $\rat_v = \rat_{p_v}$.

Let $\|\cdot \|_v$ be the unique extension to $K_v$ of the usual $p_v$-adic absolute value on $\rat_v$, and given a point $\alpha\in K$, we define
\begin{equation} \label{GNorm}
	\|\alpha\|_{\mathcal G} = \sum_{v\in M_K} \frac{[K_v:\rat_v]}{[K:\rat]} \cdot \big| \log \|\alpha\|_v\big|.
\end{equation}
The right hand side of \eqref{GNorm} does not depend on $K$, and moreover, its value is unchanged when $\alpha$ is multiplied by a root of unity.  
The properties of absolute values now imply that $\|\cdot\|_{\mathcal G}$ is a norm on $\mathcal G$ with respect to the usual absolute value on $\rat$.  By the product formula, the map $\alpha\mapsto \frac{1}{2} \|\alpha\|_{\mathcal G}$
is equal to the Weil height as defined in \cite{AllcockVaaler}.

As part of their breakthrough article \cite{AllcockVaaler}, Allcock and Vaaler showed how to represent the completion of $\mathcal G$ as a certain subspace of an $L^1$ space.  Specifically, let $Y$ denote the set of 
all places of $\alg$, and for each number field $K$ and each place $v$ of $K$, write
\begin{equation*}
	Y(K,v) = \left\{ y\in Y: y\mid v\right\}.
\end{equation*}
For every place $y\in Y(K,v)$, there exists an absolute value $\|\cdot\|_y$ on $\alg$ that extends the absolute value $\|\cdot\|_v$ on $K$.  The authors of \cite{AllcockVaaler} defined a locally compact, totally disconnected, Hausdorff topology 
on $Y$ having the set $\{ Y(K,v): [K:\rat] < \infty,\ v\in M_K\}$ as a basis.  They further defined a regular measure $\lambda$ on the Borel sets $\mathcal B$ of $Y$ with the property that
\begin{equation} \label{LambdaMeasure}
	\lambda(Y(K,v)) = \frac{[K_v:\rat_v]}{[K:\rat]}
\end{equation}
for each number field $K$ and each place $v\in M_K$.  If $\alpha\in \mathcal G$ then we let $f_\alpha:Y\to \real$ be given by 
\begin{equation*}
	f_\alpha(y) = \log \|\alpha\|_y,
\end{equation*}
and therefore, the value $\|\alpha\|_{\mathcal G}$ may be rewritten as
\begin{equation*}
	\|\alpha\|_{\mathcal G} = \int_Y |f_\alpha(y)| d\lambda(y).
\end{equation*}
Gubler \cite{Gubler} achieved something similar several years prior to \cite{AllcockVaaler}, expressing certain height functions as integrals.
Allcock and Vaaler proved that $\alpha\mapsto f_\alpha$ is an isometric isomorphism of $\mathcal G$ onto a dense $\rat$-linear subspace of
\begin{equation} \label{Y_SDefinition}
	\left\{ f\in L^1(Y,\mathcal B,\lambda): \int_Y f(y)d\lambda(y) = 0\right\}.
\end{equation}
Over the last several years, a variety of authors have applied the methods of \cite{AllcockVaaler} to establish new results on $\mathcal G$ and related structures
(see \cite{AkhtariVaaler,FiliMinerNorms,FiliMinerDirichlet,FiliMinerOrthogonal,GrizzardVaaler,Vaaler,Hughes,KellySamuels}, for example).

The author \cite{Samuels} recently began explorations into various dual spaces related to $\mathcal G$.  Let $Y^0$ be the set of non-Archimedean places of $\alg$ and define
\begin{equation*}
	\overline \intg = \left\{ \alpha\in \alg: \|\alpha\|_y \leq 1\mbox{ for all } y\in Y^0\right\}.
\end{equation*}
By the strong triangle inequality, $\overline \intg$ is a subring of $\alg$, and the group of units of $\overline \intg$ is equal to 
\begin{equation*}
	{\overline \intg}^\times =  \left\{ \alpha\in \algt: \|\alpha\|_y  = 1\mbox{ for all } y\in Y^0\right\}.
\end{equation*}
The set $\mathcal V := \algt/{\overline \intg}^\times$ is a vector space over $\rat$ under the operations analogous to \eqref{GOperations}.

For a number field $K$, we write $M^0_K$ to denote the set of non-Archimedean places of $K$.
If $\alpha\in K$ then we define
\begin{equation} \label{VNorm}
	\|\alpha\|_{\mathcal V} = \sum_{v\in M^0_K} \frac{[K_v:\rat_v]}{[K:\rat]} \cdot  \big| \log \|\alpha\|_v\big|.
\end{equation}
Analogous to $\mathcal G$, we obtain that $\|\cdot\|_{\mathcal V}$ is a norm on $\mathcal V$ with respect to the usual absolute value on $\rat$.  

The main results of \cite{Samuels} established representation theorems for the algebraic and continuous duals of $\mathcal V$ (see \cite[Theorems 1.1 \& 1.2]{Samuels}).  The primary tool used in both proofs is an object called a consistent map
which we shall briefly describe here.  Let
\begin{equation*}
	\mathcal I = \left\{(K,v): [K:\rat] < \infty,\ v\in M_K^0\right\}.
\end{equation*}
A map $c:\mathcal I\to \rat$ is called {\it consistent} if we have 
\begin{equation*} \label{ConsistentDefinitionOld}
	c(K,v) = \sum_{w\mid v} c(L,w)
\end{equation*}
for all number fields $K$, all non-Archimedean places $v$ of $K$, and all finite extensions $L/K$.   If $c,d:\mathcal I\to \rat$ are consistent maps and $r\in \rat$ then we let
\begin{equation} \label{ConsistentArithmetic}
	(c+d)(K,v) = c(K,v) + d(K,v)\quad\mbox{and}\quad (rc)(K,v) = rc(K,v).
\end{equation}
It is easily verified that these operations cause the set of consistent maps $c:\mathcal I\to \rat$ to be a vector space over $\rat$, which we shall denote by $\mathcal I^*$.  Additionally,
we defined the subspace $\mathcal I'$ of $\mathcal I^*$ by
\begin{equation} \label{ClassicBounded}
	\mathcal I' = \left\{ c\in \mathcal I^*: \frac{[K:\rat]\cdot c(K,v)}{[K_v:\rat_v]\cdot \log p_v}\mbox{ is bounded for } (K,v)\in \mathcal I\right\}.
\end{equation}
The results of \cite{Samuels} identified a particular vector space isomorphism from $\mathcal I^*$ to the algebraic dual $\mathcal V^*$ of $\mathcal V$.
Moreover, we showed that the restriction of this map to $\mathcal I'$ is itself an isomorphism onto the continuous dual $\mathcal V'$ of $\mathcal V$.  

With some appropriate minor adjustments, consistency turns out to be precisely the property required to classify duals of several locally constant function spaces on places of $\alg$ arising in \cite{AllcockVaaler}.  
The purpose of the present article is to provide the relevant family of representation theorems, yielding four main results (Theorems \ref{MainIsomorphism}, \ref{ZeroIsomorphism}, \ref{ContinuousIsomorphism} and \ref{ContinuousZeroIsomorphism}).  
We regard these theorems as versions of the Riesz Representation Theorem analogous to those found in \cite{Rudin,DunfordSchwartz,Bogachev}, but with consistent maps playing the role of measures.  
Additionally, our work generalizes \cite[Theorem 1.2]{Samuels}.

\subsection{Main Results} 
Let $S$ be a set of places of $\rat$ and $Y_S$ the set of places of $\alg$ which divide a place in $S$.  Extending our definition of $\mathcal I$ given above, we let
\begin{equation*}
	\mathcal J_S = \left\{(K,v): [K:\rat] < \infty,\ v\in M_{K,S}\right\},
\end{equation*}
where $M_{K,S}$ denotes the set of places of $K$ which divide a place in $S$.  A map $c:\mathcal J_S\to \real$ is called {\it consistent} if it satisfies 
\begin{equation} \label{ConsistentDefinition}
	c(K,v) = \sum_{w\mid v} c(L,w)
\end{equation}
for all number fields $K$, all places $v\in M_{K,S}$, and all finite extensions $L/K$.  Note that our new definition of consistent permits $c$ to take real values, as opposed to our previous definition, which considered rational valued functions only.
With addition and scalar multiplication defined as in \eqref{ConsistentArithmetic}, the collection of consistent maps forms an $\real$-vector space which we shall denote by $\mathcal J_S^*$.

The collection $\mathfrak B_S := \{Y(K,v): (K,v)\in \mathcal J_S\}$ forms a basis for the topology on $Y_S$ given in \cite{AllcockVaaler}.  If $\mu$ is a signed Borel measure on $Y_S$, then the elementary
properties of measures ensure that $(K,v)\mapsto \mu(Y(K,v))$ defines a consistent map.   However, our definition of consistent map requires only finite additivity on $\mathfrak B_S$ rather than the countable 
additivity that is required of measures.  Therefore, not every consistent map arises from a signed Borel measure in this manner.

One particularly useful consistent map arises from the measure $\lambda$ on $Y_S$ defined in \cite{AllcockVaaler}.  We shall simply write $\lambda:\mathcal J_S\to \real$ to denote the associated consistent map so that
\begin{equation*}
	\lambda(K,v) = \frac{[K_v:\rat_v]}{[K:\rat]}.
\end{equation*}

%We now write $\lcc$ to denote the $\real$-vector space of locally constant functions from $Y_S$ to $\real$ having compact support.  We advise the reader that $LC_c(Y_S)$ is often used in other parts of the literature
%to denote this space (see \cite{AllcockVaaler,Rudin}, for example).  In those settings, such a robust notation is necessary because this space appears alongside other function spaces such as $C_c(Y_S)$.
%Our results, on the other hand, pertain only to spaces of locally constant functions, and therefore, we prefer the simplified notation $\lcc$.  
Following the notation of \cite{Rudin}, we write $\lcc$ to denote $\real$-vector space of locally constant functions from $Y_S$ to $\real$ having compact support.
In the special case where $Y_S= Y$, the results of \cite[\S 4]{AllcockVaaler} contain 
several important results on $\lcc$ that are needed to establish their main results.   In particular, we note Lemma 5, which establishes that $\lcc$ is dense in $L^p(Y,\mathcal B,\lambda)$ for every $1\leq p < \infty$.  
The following proposition regarding elements of $\lcc$ follows essentially from \cite[Lemma 4]{AllcockVaaler}.
 
\begin{prop}\label{CompactnessRestriction}
	If $f\in \lcc$ then there exists a number field $K$ such that $f$ is constant on $Y(K,v)$ for all $v\in M_{K,S}$.  Moreover, that constant is equal to $0$ for all but finitely many $v\in M_{K,S}$.
\end{prop}

We shall write $\Omega_S(f)$ to denote the set of number fields which satisfy the conclusions of Proposition \ref{CompactnessRestriction}.  If $K\in\Omega_S(f)$ and $L$ is a finite extension of $K$,
then we obtain easily that $L\in \Omega_S(f)$.  Moreover, if $f\in \lcc$ and $K\in \Omega_S(f)$, we may interpret $f$ as a well-defined function $f:M_{K,S}\to \real$ satisfying the following two properties:
\begin{enumerate}[label={(\roman*)}]
	\item $f(v) = f(y)$ for all $y\in Y(K,v)$
	\item $f(v) = 0$ for all but finitely many $v\in M_{K,S}$.
\end{enumerate}
Given a consistent map $c\in \mathcal J_S^*$ and $K\in \Omega_S(f)$, we let
\begin{equation} \label{Phic}
	\Phi_c(f) = \sum_{v\in M_{K,S}} f(v)c(K,v).
\end{equation}
The definition of consistency implies that the right hand side of \eqref{Phic} is independent of the choice of $K\in \Omega_S(f)$.  As $\Phi_c$ is easily verified to be linear, it is a well-defined
element of the algebraic dual $\lcc^*$ of $\lcc$.  We let $\Phi_S^*:\mathcal J_S^*\to \lcc^*$ be given by $\Phi_S^*(c) = \Phi_c$ and state our first main result -- a representation theorem for $\lcc^*$ in terms of consistent maps.

\begin{thm} \label{MainIsomorphism}
	$\Phi_S^*$ is a vector space isomorphism from $\mathcal J_S^*$ to $\lcc^*$.
\end{thm}

Given a point $\alpha\in \alg$, we recall that Allcock and Vaaler \cite{AllcockVaaler} defined $f_\alpha:Y\to \real$ by $f_\alpha(y) = \log \|\alpha\|_y$.
From this point forward, we shall write $f_{S,\alpha}$ to denote the restriction of this map to $Y_S$, i.e., $f_{S,\alpha}:Y_S\to \real$ is given by 
\begin{equation*}
	f_{S,\alpha}(x) = \log\|\alpha\|_x.
\end{equation*}  
We clearly have that $f_{S,\alpha}\in \lcc$, however, when $S$ is the set of all places of $\rat$, something more is true.
Indeed, the product formula implies that $f_{S,\alpha}$ belongs to the codimension $1$ subspace of $\lcc$ given by
\begin{equation*}
	\lcco := \left\{ f\in \lcc: \int_{Y_S} f(x)d\lambda(x) = 0\right\}.
\end{equation*}
As a result, we would like a representation theorem for $\lcco^*$ which is analogous to Theorem \ref{MainIsomorphism}.

If $c:\mathcal J_S\to \real$ is a consistent map, then $\Phi_c$ may be restricted to $\lcco$ to obtain a well-defined element of $\lcco^*$.  Therefore, we may let $\Psi_S^*:\mathcal J_S^*\to \lcco^*$ 
be given by 
\begin{equation*}
	\Psi_S^*(c) = \Phi_c\big|_{\lcco}.
\end{equation*}
For a given consistent map $c$, it is possible that $\Psi_S^*(c)$ is trivial even if $\Phi_S^*(c)$ is not.  Indeed, we could have that $\Phi_c(f) = 0$ for all $f\in \lcco$ but $\Phi_c(g) \ne 0$ for some $g\in \lcc$.
Our next result shows precisely when this phenomenon occurs.  We remind the reader that $\lambda\in \mathcal J_S^*$ is the particular consistent map given by $\lambda(K,v) = [K_v:\rat_v]/[K:\rat]$.

\begin{thm} \label{ZeroIsomorphism}
	$\Psi_S^*$ is a surjective linear transformation from $\mathcal J_S^*$ to $\lcco^*$ such that $\ker\Psi_S^* = \spann\{\lambda\}$.
\end{thm}

We now let $\mathcal F_S = \{f_{S,\alpha}:\alpha\in \alg\}$ and note that $\mathcal F_S$ is a vector space over $\rat$ and a subset of $\lcc$.  This space is vaguely familiar from \cite{AllcockVaaler} and \cite{Samuels}.
Indeed, if $S$ is the set of all places of $\rat$ then $\alpha\mapsto f_{S,\alpha}$ is an isomorphism of $\mathcal G$ with $\mathcal F_S$.  Similarly, if $S$ is the set of non-Archimedean places
of $\rat$ then $\alpha\mapsto f_{S,\alpha}$ is an isomorphism of $\mathcal V$ with $\mathcal F_S$.

Let $\|\cdot \|_1$ denote the $L^1$ norm on $\lcc$ so that
\begin{equation} \label{L1Norm}
	\|f\|_1 = \int_{Y_S} |f(x)|d\lambda(x).
\end{equation}
If $K\in \Omega_S(f)$, where $\Omega_S(f)$ is defined in the remarks following Proposition \ref{CompactnessRestriction}, then we may write $f(v) = f(y)$ for all $y\in Y(K,v)$.  In this scenario, we obtain an alternate expression for $\|f\|_1$ given by
\begin{equation} \label{AltNorm}
	\|f\|_1 = \sum_{v\in M_{K,S}} \frac{[K_v:\rat_p]}{[K:\rat]} |f(v)|.
\end{equation}
The following observations may now be obtained from \cite{AllcockVaaler}.

\begin{prop} \label{FSubspace}
	$\mathcal F_S$ is a $\rat$-linear subspace of $\lcc$ satisfying the following properties:
	\begin{enumerate}[label={(\roman*)}]
		\item\label{AllY} If $Y_S=Y$ then $\mathcal F_S$ is dense in $\lcco$.
		\item\label{SomeY} If $Y_S\subsetneq Y$ then $\mathcal F_S$ is dense in $\lcc$.
	\end{enumerate}
\end{prop}

The article \cite{Samuels} examined the algebraic dual $\mathcal F_S^*$ of $\mathcal F_S$ in the case where $S$ is the set of non-Archimedean places of $\rat$.  Specifically, \cite[Theorem 1.1]{Samuels} provides a representation theorem
for $\mathcal F_S^*$ that is analogous to Theorem \ref{MainIsomorphism}.  It is unlikely, however, that Theorems \ref{MainIsomorphism} and \ref{ZeroIsomorphism} can be combined with Proposition \ref{FSubspace} to yield a new proof of 
\cite[Theorem 1.1]{Samuels}.  Indeed, all of the aforementioned results deal exclusively with algebraic dual spaces, and hence, it is difficult to imagine capitalizing on the density properties of Proposition \ref{FSubspace} to say something about
$\mathcal F_S^*$.  In view of these observations, we find it useful to examine the continuous dual spaces $\lcc'$ and $\lcco'$ of $\lcc$ and $\lcco$, respectively.  

To this end, we define
\begin{equation*}
	\mathcal J_S' = \left\{ c\in \mathcal J_S^*: \frac{c(K,v)}{\lambda(K,v)}\mbox{ is bounded for } (K,v)\in \mathcal J_S\right\}
\end{equation*}
and note that $\mathcal J_S'$ is a subspace of $\mathcal J_S^*$ containing $\lambda$.  Let $\Phi_S'$ denote the restriction of $\Phi_S^*$ to $\mathcal J_S'$.  Our next result is a continuous analog of Theorem \ref{MainIsomorphism}.

\begin{thm} \label{ContinuousIsomorphism}
	$\Phi_S'$ is a vector space isomorphism from $\mathcal J_S'$ to $\lcc'$.
\end{thm}

If $\Phi\in \lcc'$ then clearly the restriction of $\Phi$ to $\lcco$ is also continuous, and hence, we may define
$\Psi_S':\mathcal J_S'\to \lcco'$ by 
\begin{equation*}
	\Psi_S'(c) = \Phi_c\big|_{\lcco}.
\end{equation*}  
Equivalently, we may interpret $\Psi_S'(c)$ as the restriction of $\Phi_S'(c)$ to $\lcco$.  Not surprisingly, we obtain a continuous analog of Theorem \ref{ZeroIsomorphism}.

\begin{thm} \label{ContinuousZeroIsomorphism}
	$\Psi_S'$ is a surjective linear transformation from $\mathcal J_S'$ to $\lcco'$ such that $\ker\Psi_S' = \spann\{\lambda\}$.
\end{thm}

The remainder of this paper is organized in the following way.  In Sections \ref{Algebraic} and \ref{Continuous}, we present the proofs of our main results, with Theorems \ref{MainIsomorphism} and \ref{ZeroIsomorphism}
coming in Section \ref{Algebraic}.  We use Section \ref{Applications} to discuss how our main results relate to those of \cite{Samuels} and also to pose several problems left open by our work.  
Finally, the reader may have noticed that, in broadening the definition of consistent from \cite{Samuels}, we could have taken $S$ to be a set of places of an arbitrary number field $F$ rather than restricting our attention to $\rat$.  
In Section \ref{Generalization}, we posit this alternate definition of consistency, and subsequently, we show it to be essentially equivalent to the definition we have already provided.

\section{Algebraic Dual Spaces: Proofs of Theorems \ref{MainIsomorphism} and \ref{ZeroIsomorphism}} \label{Algebraic}

As we noted earlier, Proposition \ref{CompactnessRestriction} is a fairly standard compactness proof that essentially follows from \cite[Lemma 4]{AllcockVaaler}.  Even though the differences are only trivial, we include its proof here for the
sake of completeness.  As was the case in our introduction, we continue to assume that $S$ is a fixed set of places of $\rat$ and $Y_S$ is the set of places of $\alg$ which divide a place in $S$.

\begin{proof}[Proof of Proposition \ref{CompactnessRestriction}]
	Clearly the set $\{Y(\rat,p): p\in S\}$ is a cover of $\supp(f)$, so by compactness, there exists a finite subset $T\subseteq S$ such that
	\begin{equation} \label{SupportContain}
		\supp(f) \subseteq \bigcup_{p\in T} Y(\rat,p).
	\end{equation}
	Letting $W$ denote the right hand side of \eqref{SupportContain}, we note that $W$ must be compact since it is a finite union of compact sets.
	
	We have assumed that $f$ is locally constant, so for each $y\in W$ there exists an open neighborhood $U_y$ of $Y$ such that $f$ is constant on $U_y$.  By definition of the topology on $Y$, we may assume that
	$U_y = Y(K_y,v_y)$ for some number field $K_y$ and some place $v_y$ of $K_y$.  We clearly have that $y\mid v_y$.  Certainly $\{U_y:y\in W\}$ is a cover of $W$, and since $W$ is compact, there
	exists a finite subset $W_0\subseteq W$ such that
	\begin{equation}\label{OpenCover}
		\{U_y:y\in W_0\} 
	\end{equation}
	is a cover of $W$.  Now let $K$ be the compositum of $K_y$ for $y\in W_0$. We claim that $f$ is constant on $Y(K,v)$ for all $v\in M_{K,S}$.
	
	If $v$ does not divide a place in $T$, then $Y(K,v) \subseteq Y_S\setminus W$.  In this case, $f$ is equal to $0$ on $Y(K,v)$ so that the required conclusion holds.
	Since there are only finitely many places of $K$ dividing a place in $T$, the second statement of the proposition follows from this observation as well.
	
	Now suppose $v$ divides a place $p\in T$ and let $x\in Y(K,v)$.  Since \eqref{OpenCover} is a cover of $W$, there must exist $y\in W_0$ such that 
	\begin{equation*}
		x \in U_y = Y(K_y,v_y).
	\end{equation*}
	Since $K$ is an extension of $K_y$, we have now shown that $v\mid v_y$, and so $Y(K,v)\subseteq Y(K_y,v_y)$.  It now follows that $f$ is constant on $Y(K,v)$.
\end{proof}

We recall that $\Omega_S(f)$ denotes the set of number fields which satisfy the conclusions of Proposition \ref{CompactnessRestriction}, and if $f\in \lcc$ and $K\in \Omega_S(f)$, we may regard $f$ as a well-defined 
function $f:M_{K,S}\to \real$ satisfying the following:
\begin{enumerate}[label={(\roman*)}]
	\item $f(v) = f(y)$ for all $y\in Y(K,v)$
	\item $f(v) = 0$ for all but finitely many $v\in M_{K,S}$.
\end{enumerate}
For any map $c:\mathcal J_S\to \real$ we may now define
\begin{equation} \label{PhiDefInitial}
	\Phi_c(K,f) = \sum_{v\in M_{K,S}} f(v)c(K,v).
\end{equation}
The reader may recognize the right hand side of \eqref{PhiDefInitial} from our definition of $\Phi_c(f)$ in \eqref{Phic}.  However, as we have not assumed $c$ to be consistent, we may not assume that \eqref{PhiDefInitial} is independent of $K$.
Whether or not $c$ is consistent plays an important role in studying this definition.

\begin{lem} \label{ConsistencyIndependence}
	Suppose that $c:\mathcal J_S\to \real$ is any map.  Then $c$ is consistent if and only if $\Phi_c(K,f) =\Phi_c(L,f)$ for every $f\in \lcc$ and every $K,L\in \Omega_S(f)$.
\end{lem}	
\begin{proof}
	We first suppose that $c$ is consistent.  Suppose that $f\in \lcc$ and $K,L\in \Omega_S(f)$.  By possibly replacing $L$ by the compositum of $K$ and $L$, it is sufficient to consider the case where $K\subseteq L$.  
	If $w$ is a place of $L$ dividing the place $v$ of $K$, then $f(w) = f(v)$.  We now obtain
	\begin{align*}
		\sum_{w\in M_L(S)} c(L,w)f(w) & = \sum_{v\in M_{K,S}} \sum_{w\mid v} c(L,w)f(w) \\
								& = \sum_{v\in M_{K,S}} f(v) \sum_{w\mid v} c(L,w) \\
								& =  \sum_{v\in M_{K,S}} c(K,v) f(v)
	\end{align*}
	as required.
	
	Next assume that $\Phi_c(K,f) =\Phi_c(L,f)$ for every $f\in \lcc$ and every $K,L\in \Omega_S(f)$.  In order to show that $c$ is consistent, we let $(K,v)\in \mathcal J_S$ and $L$ a finite extension of $K$.
	Also let $f\in \lcc$ be the indicator function of $Y(K,v)$, i.e., 
	\begin{equation*}
		f(y) = \begin{cases} 1 & \mbox{if } y\mid v \\ 0 & \mbox{if } y\nmid v.\end{cases}
	\end{equation*}
	Clearly $K$ and $L$ both belong to $\Omega_S(f)$, so our assumption implies that $\Phi_c(K,f) = \Phi_c(L,f)$.  This observation yields
	\begin{align*}
		c(K,v) =  \Phi_c(L,f) & = \sum_{w\in M_L(S)} c(L,w)f(w) \\
						& = \sum_{v\in M_{K,S}} f(v) \sum_{w\mid v} c(L,w) \\
						& = \sum_{w\mid v} c(L,w)
	\end{align*}
	proving that $c$ is consistent.
\end{proof}

We are now prepared to present our proofs of Theorems \ref{MainIsomorphism} and \ref{ZeroIsomorphism}.

\begin{proof}[Proof of Theorem \ref{MainIsomorphism}]
	It is straightforward to check that $\Phi_S^*$ is linear, so we begin by assuming $c\in \mathcal J_S^*$ is such that $\Phi_c \equiv 0$.  To see that $c\equiv 0$, we let $(K,v)\in \mathcal J_S$.  Now define 
	\begin{equation} \label{Indicator}
		f_v(y) = \begin{cases} 1 & \mbox{if } y\mid v \\ 0 & \mbox{if } y\nmid v \end{cases}
	\end{equation}
	so that $f_v$ is a well-defined element of $\lcc$.  Since we have assumed $\Phi_c \equiv 0$, we conclude that
	\begin{equation*}
		0 = \sum_{u\in M_{K,S}} c(K,u)f_v(u) = c(K,v)
	\end{equation*}
	which shows that $c\equiv 0$.  We have now shown that $\Phi_S^*$ is injective.
	
	We now proceed with the proof that $\Phi_S^*$ is surjective.  We assume that $\Phi\in \lcc^*$ and seek an element $c\in \mathcal J_S^*$ such that $\Phi_c =\Phi$.
	For each number field $K$, we let
	\begin{equation*}
		T_K = \left\{(a_v)_{v\in M_{K,S}}: a_v\in \real,\ a_v =0 \mbox{ for all but finitely many } v\in M_{K,S}\right\}
	\end{equation*}
	and note that $T_K$ is clearly a vector space over $\real$.  Further let $\zz_v$ be the element of $T_K$ having a $1$ at entry $v$ and $0$ elsewhere.  Clearly $\{\zz_v:v\in M_{K,S}\}$ is a 
	basis for $T_K$ over $\real$.  Also let $\mathcal C_K$ be the subset of $\lcc$ consisting of those functions which are constant on $Y(K,v)$ for all $v\in M_{K,S}$.  If $f\in \mathcal C_K$ and $y\in Y(K,v)$ we shall write $f(v) = f(y)$. 
	Clearly $\mathcal C_K$ is a subspace of $\lcc$.
	
	We define $\Delta_K:\mathcal C_K\to T_K$ by 
	\begin{equation*}
		\Delta_K(f) = (f(v))_{v\in M_{K,S}}
	\end{equation*}
	and claim that $\Delta_K$ is an isomorphism.  Clearly $\Delta_K$ is linear.  Moreover, if $\Delta_K(f) = 0$ then $f(v) = 0 $ for all $v\in M_{K,S}$.  But $f$ is constant on $Y(K,v)$ for all $v\in M_{K,S}$, which implies
	that $f(y) = 0$ for all $y\in Y_S$.  This shows that $\Delta_K$ is injective.  Certainly $f_v\in \mathcal C_K$ and $\Delta_K(f_v) = \zz_v$, which shows that $\zz_v$ is in the range of $\Delta_K$.  As $\{\zz_v\}$ forms a basis for $T_K$ over $\real$,
	it follows that $\Delta_K$ is surjective.
	
	Now let $(K,v)\in \mathcal J_S$.  Let $A_K: T_K\to \real$ be given by $A_K = \Phi\circ \Delta_K^{-1}$ and define
	\begin{equation*}
		c(K,v) = A_K(\zz_v).
	\end{equation*}
	To complete the proof, we show that $c$ is consistent and $\Phi_c = \Phi$.
	
	For an arbitrary point $\xx = (x_v)_{v\in M_{K,S}} \in T_K$, we have 
	\begin{equation} \label{LKFormula}
		A_K(\xx) = \sum_{v\in M_{K,S}} x_v A_K(\zz_v) = \sum_{v\in M_{K,S}} x_vc(K,v).
	\end{equation}
	Given an element $f\in \lcc$, we apply Lemma \ref{CompactnessRestriction} to select a number field $K\in \Omega_S(f)$.  This means that $f\in \mathcal C_K$ and we may apply \eqref{LKFormula} with $\xx = \Delta_K(f)$ to obtain
	\begin{equation} \label{PhiRewrite}
		\Phi(f) = A_K(\Delta_K(f)) = \sum_{v\in M_{K,S}} f(v) c(K,v) = \Phi_c(K,f).
	\end{equation}
	Of course, if $L$ is a different number field in $\Omega_S(f)$ then the same argument shows that $\Phi(f) = \Phi_c(L,f)$.  Now Lemma \ref{ConsistencyIndependence} means that $c$ is consistent and $\Phi_c(f) = \Phi(f)$
	as required.
\end{proof}

\begin{proof}[Proof of Theorem \ref{ZeroIsomorphism}]
	Suppose $\iota^*:\lcc^*\to \lcco^*$ is the restriction map, i.e., $\iota^*(f)$ is equal to the restriction of $f$ to $\lcco$.  If we wish, we may interpret $\iota^*$ as the pullback of the inclusion map 
	$\iota:\lcco\to \lcc$ given by $\iota^*(f) = f\circ \iota$.  Clearly $\iota^*$ is a linear map and
	\begin{equation*}
		\Psi_S^* = \iota^*\circ \Phi_S^*.
	\end{equation*}
	Since $\Phi_S^*$ is an isomorphism (see Theorem \ref{MainIsomorphism}) and $\iota^*$ is linear, it follows that $\Psi_S^*$ is linear.
	
	In order to prove that $\Psi_S^*$ is surjective, it is sufficient to show that $\iota^*$ is surjective.  We first show that there exists a one-dimensional subspace $\mathcal W\subseteq \lcc$ such that 
	\begin{equation} \label{DirectSum}
		\lcc = \lcco \oplus \mathcal W.
	\end{equation}
	To see this, we select a place $p\in S$, and for each $r\in \real$, we define
	\begin{equation*}
		f_r(y) = \begin{cases} r & \mbox{if } y\in Y(\rat,p) \\ 0 & \mbox{if } y\not\in Y(\rat,p).\end{cases}
	\end{equation*}
	Now set $\mathcal W = \{f_r:r \in \real\}$ so that $\mathcal W$ is clearly a one-dimensional subspace of $\lcc$ and $\mathcal W\cap\lcco = \{0\}$.  It remains only to show that $\lcc = \lcco + \mathcal W$.
	Let $g\in \lcc$ and let
	\begin{equation*}
		r = \int_{Y_S} g(x)d\lambda(x).
	\end{equation*}
	Then define $h(y) = g(y) - f_r(y)$ and observe that
	\begin{align*}
		\int_{Y_S} h(y)d\lambda(y) & = \sum_{q\in S} \int_{Y(\rat,q)} h(y)d\lambda(y) \\
			& = \int_{Y(\rat,p)} (g(y) - r)d\lambda(y) + \sum_{q\in S\setminus \{p\}} \int_{Y(\rat,q)} g(y)d\lambda(y) \\
			& = \int_{Y_S} g(y)d\lambda(y) - \int_{Y(\rat,p)} r d\lambda(y)  = 0,
	\end{align*}
	where the last equality follows from the definition of $r$.  We have now shown that $g = h + f_r\in \lcco + \mathcal W$ and \eqref{DirectSum} follows.
	Based on this observation, each element $g\in \lcc$ may be uniquely represented in the form $g = h + f_r$, where $h\in \lcco$ and $r\in \real$.
	
	As a result, if $\phi\in \lcco^*$ we may define $\Phi:\lcc\to \real$ by $\Phi(g) = \phi(h)$.  Certainly $\Phi$ is linear and $\Phi(h) = \phi(h)$ for all $h\in \lcco$.
	For each such element $h$, we have 
	\begin{equation*}
		(\iota^*(\Phi))(h) = (\Phi\circ\iota)(h) = \Phi(h) = \phi(h)
	\end{equation*}
	which shows that $\iota^*(\Phi) = \phi$, establishing that $\iota^*$ is surjective.
	
	We now proceed with our proof that $\ker\Psi_S^* = \spann\{\lambda\}$.  Assume that $f\in \lcco$ and suppose that $K\in \Omega_S(f)$ so we have
	\begin{equation*}
		\Phi_\lambda(f) = \sum_{v\in M_{K,S}} \lambda(K,v)f(v) = \sum_{v\in M_{K,S}} \frac{[K_v:\rat_v]}{[K:\rat]} f(v).
	\end{equation*}
	Using the definition of the measure $\lambda$ on $Y_S$ along with the fact that $f$ is constant on $Y(K,v)$, we obtain
	\begin{equation*}
		\int_{Y(K,v)} f(x)d\lambda(x) = f(v)\lambda(Y(K,v)) = \frac{[K_v:\rat_v]}{[K:\rat]} f(v).
	\end{equation*}
	It now follows that
	\begin{equation*}
		\Phi_\lambda(f) =  \sum_{v\in M_{K,S}} \int_{Y(K,v)} f(x)d\lambda(x) = \int_{Y_S} f(x)d\lambda(x) = 0.
	\end{equation*}
	We have shown that $\Psi_S^*(\lambda) = 0$ implying that $\spann\{\lambda\} \subseteq \ker\Psi_S^*$.
	
	Now suppose that $c\in \ker\Psi_S^*$ which means that $\Phi_c(h) = 0$ for all $h\in \lcco$.  For each number field $K$, we continue to write $\mathcal C_K$ for the subspace of $\lcc$ consisting of those
	functions that are constant on $Y(K,v)$ for all $v\in M_{K,S}$.  If $h\in \mathcal C_K\cap \lcco$ then our assumptions ensure that
	\begin{equation} \label{ZeroH}
		\sum_{v\in M_{K,S}} c(K,v)h(v) = 0.
	\end{equation}
	We claim that for every number field $K$, there exists $r_K\in \real$ such that
	\begin{equation} \label{LocalMultiple}
	 	c(K,v) = r_K \cdot\frac{[K_v:\rat_v]}{[K:\rat]}\quad\mbox{for all }v\in M_{K,S}.
	\end{equation}
	If $M_{K,S}$ contains only one place $u$, then \eqref{LocalMultiple} follows easily by setting $r_K = [K:\rat] c(K,u) / [K_u:\rat_u]$.  We assume that $M_{K,S}$ contains at least two places
	and fix a place $u\in M_{K,S}$.  Now setting $r_u = [K:\rat] c(K,u) / [K_u:\rat_u]$, we clearly have that 
	\begin{equation*}
	 	c(K,u) = r_u \cdot\frac{[K_u:\rat_u]}{[K:\rat]}.
	\end{equation*}
	Now let $v\in M_{K,S}$ be such that $v\ne u$ and define 
	\begin{equation*}
		h(y) = \begin{cases} [K:\rat]/[K_u:\rat_u] & \mbox{if } y\mid u \\ 
						-[K:\rat]/[K_v:\rat_v] & \mbox{if } y\mid v \\
						0 & \mbox{otherwise}.  \end{cases}
	\end{equation*}
	One easily checks that $h\in \mathcal C_K\cap \lcco$, so \eqref{ZeroH} applies to yield
	\begin{equation*}
		0 = c(K,u)h(u) + c(K,v)h(v) = r_u \frac{[K_u:\rat_u]}{[K:\rat]}\cdot \frac{[K:\rat]}{[K_u:\rat_u]} - c(K,v) \frac{[K:\rat]}{[K_v:\rat_v]}.
	\end{equation*}
	We now obtain that
	\begin{equation*}
		c(K,v) = r_u \cdot\frac{[K_v:\rat_v]}{[K:\rat]} \quad\mbox{for all }v\in M_{K,S}.
	\end{equation*}
	For each number field $K$, we have now found $r_K\in \real$ such that \eqref{LocalMultiple} holds.
	
	To complete the proof, we claim that $r_K$ is independent of $K$.  To see this, observe that $r_\rat = c(\rat,p)$ for all $p\in S$.  If $K$ is any number field, then the consistency of $c$ implies that
	\begin{equation*}
		r_\rat = \sum_{v\mid p} c(K,v) = \sum_{v\mid p} r_K \cdot\frac{[K_v:\rat_v]}{[K:\rat]} = r_K.
	\end{equation*}
	In other words
	\begin{equation*} 
	 	c(K,v) = r_\rat\cdot\frac{[K_v:\rat_v]}{[K:\rat]}\quad\mbox{for all }(K,v)\in \mathcal J_S,
	\end{equation*}
	which means that $c\in \spann\{\lambda\}$ as required.
\end{proof}

\section{Continuous Dual Spaces: Proofs of Theorems \ref{ContinuousIsomorphism} and \ref{ContinuousZeroIsomorphism}} \label{Continuous}

Although the proof of Proposition \ref{FSubspace} is not required to prove Theorems \ref{ContinuousIsomorphism} and \ref{ContinuousZeroIsomorphism}, we include it here as it relates most closely to those results.

\begin{proof}[Proof of Proposition \ref{FSubspace}]
	Define the $\rat$-vector space $\mathcal G = \alg/\algt_\tors$ and note that $f_{S,\alpha}(y)$ is well-defined for $\alpha\in \mathcal G$.  Now define $\phi:\mathcal G\to \lcc$ by $\phi(\alpha) = f_{S,\alpha}$ so that $\phi(\mathcal G) = \mathcal F_S$.
	One easily checks that $\phi$ is a $\rat$-linear map, so the first statement of the proposition follows immediately, and then \ref{AllY} follows from \cite[Theorem 1]{AllcockVaaler}.
	
	To prove \ref{SomeY}, we may assume that $p$ is a place of $\rat$ not belonging to $S$.  Now let $g\in \lcc$ and $\varepsilon > 0$.  Also set
	\begin{equation*}
		C = \int_{Y_S} g(x)d\lambda(x)
	\end{equation*}
	and define $h\in LC_c(Y)$ by
	\begin{equation*}
		h(x) = \begin{cases} g(x) & \mbox{if } x\in Y_S \\ -C & \mbox{if } x\mid p \\  0 & \mbox{if } x\not\in Y_S\mbox{ and } x\nmid p. \end{cases}
	\end{equation*}
	Here $LC_c(Y)$ denotes the set of locally constant functions from $Y$ to $\real$ with compact support as in \cite{Rudin}.  Using the fact that $\lambda(Y(\rat,p)) = 1$, one easily checks that
	\begin{equation*}
		\int_Y h(x)d\lambda(x) = 0.
	\end{equation*}
	By applying \ref{AllY}, there exists $\alpha\in \alg$ such that 
	\begin{equation*}
		\int_{Y_S} |f_{S,\alpha}(x) - g(x)|d\lambda(x) \leq \int_Y|f_{S,\alpha}(x) - h(x)|d\lambda(x) < \varepsilon
	\end{equation*}
	which completes the proof.
\end{proof}

We now proceed with our proofs of Theorems \ref{ContinuousIsomorphism} and \ref{ContinuousZeroIsomorphism} which rely heavily on our earlier results regarding algebraic dual spaces.

\begin{proof}[Proof of Theorem \ref{ContinuousIsomorphism}]
	From Theorem \ref{MainIsomorphism} we know that $\Phi_S^*$ is an isomorphism from $\mathcal J_S^*$ to $\lcc^*$.  In order to complete the proof, we must show that $\Phi_S^*(c)$ is continuous if and only if $c\in \mathcal J_S'$.
	Assuming that $c\in \mathcal J_S'$, there exists $B\geq 0$ such that 
	\begin{equation*}
		\frac{c(K,v)}{\lambda(K,v)} \leq B\quad\mbox{ for all }(K,v)\in \mathcal J_S.
	\end{equation*}
	Assuming that $f\in \lcc$ and $K\in \Omega_S(f)$, we have that
	\begin{equation*}
		|\Phi_c(f)| \leq \sum_{v\in M_{K,S}} |c(K,v)|\cdot |f(v)| \leq  B \sum_{v\in M_{K,S}} |\lambda(K,v)|\cdot |f(v)|  =  B\cdot  \|f\|_1
	\end{equation*}
	and it follows that $\Phi_c$ is continuous.
	
	Now assume that $c\in \mathcal J_S^*$ and $\Phi_c$ is continuous.  Since $\Phi_c$ is a continuous linear functional, there exists $B\geq 0$ such that
	\begin{equation*} 
		|\Phi_c(f)| \leq B \cdot \|f\|_1 \quad\mbox{ for all } f\in \lcc.
	\end{equation*}
	If $f\in \lcc$ and $K\in \Omega_S(f)$ then the definition of $\Phi_c$ as well as \eqref{AltNorm} implies that
	\begin{equation}\label{LipBound}
		\left| \sum_{v\in M_{K,S}} c(K,v) f(v)\right| \leq B\cdot \left | \sum_{v\in M_{K,S}} \lambda(K,v) f(v)\right|.
	\end{equation}
	Given an element $(K,v)\in \mathcal J_S$, we assume that $f_v$ is the indicator function for $Y(K,v)$.  Then we may apply \eqref{LipBound} with $f_v$ in place of $f$ to yield
	\begin{equation*}
		|c(K,v)| \leq B \cdot |\lambda(K,v)|
	\end{equation*}
	and it follows that $c\in\mathcal J_S'$ as required.
\end{proof}

\begin{proof}[Proof of Theorem \ref{ContinuousZeroIsomorphism}]
	We first observe that $\Psi_S'$ equals the restriction of $\Psi_S^*$ to $\mathcal J_S'$, so it follows from Theorem \ref{ZeroIsomorphism} that $\Psi_S'$ is linear.  Clearly $\lambda\in \mathcal J_S'$ so we conclude also that 
	$\Psi_S'(\lambda) = \Psi_S^*(\lambda) = 0$ so that $\spann\{\lambda\} \subseteq \ker\Phi_S'$.  Additionally, if $c\in \ker\Psi_S'$ then $\Psi_S^*(c) = \Psi_S'(c) = 0$.  In this case, Theorem \ref{ZeroIsomorphism} implies
	that $c\in \spann\{\lambda\}$ and we have established that $\ker\Psi_S' = \spann\{\lambda\}$.
	
	To prove that $\Psi_S'$ is surjective, we let $\phi\in \lcco'$.  In order to apply Theorem \ref{ContinuousIsomorphism}, we must prove that $\phi$ extends to a continuous linear functional on $\lcc$.
	To this end we select a place $p\in S$, and for each $r\in \real$, we define
	\begin{equation*}
		f_r(y) = \begin{cases} r & \mbox{if } y\in Y(\rat,p) \\ 0 & \mbox{if } y\not\in Y(\rat,p).\end{cases}
	\end{equation*}
	We recall from the proof of Theorem \ref{ZeroIsomorphism} that each element $g\in \lcc$ may be uniquely represented in the form $g = h + f_r$,
	where $h\in \lcco$ and
	\begin{equation*}
		r = \int_{Y_S} g(x) d\lambda(x).
	\end{equation*}
	Hence, we may define the map $\pi:\lcc\to \lcco$ by $\pi(g) = h$.  We now claim that $\pi$ is continuous with respect to the $L^1$ norm on both spaces.  We have that
	\begin{align*}
		\|\pi(g)\|_1 & = \int_{Y_S} |h(x)|d\lambda(x) \\
			& = \int_{Y_S}|h(x) - f_r(x) + f_r(x)|d\lambda(x) \\
			& \leq \int_{Y_S}|g(x)| + \int_{Y_S}|f_r(x)|d\lambda(x).
	\end{align*}
	But we have 
	\begin{equation*}
		\int_{Y_S} |f_r(x)|d\lambda(x) = |r| = \left| \int_{Y_S} g(x)d\lambda(x)\right| \leq \int_{Y_S}|g(x)|d\lambda(x)
	\end{equation*}
	and it follows that
	\begin{equation*}
		\|\pi(g)\|_1 \leq 2 \int_{Y_S}|g(x)|d\lambda(x) = 2\|g\|_1.
	\end{equation*}
	It now follows that $\pi:\lcc \to\lcco$ is continuous.  Now we define $\Phi:\lcc\to \real$ by $\Phi = \phi\circ\pi$ and we check easily that $\Phi\in \lcc'$ and that $\Phi$ is an extension of $\phi$.  
	By Theorem \ref{ContinuousIsomorphism} there exists $c\in \mathcal J_S'$ such that $\Phi_S'(c) = \Phi$, and hence, $\Psi_S'(c) = \phi$, as required.
\end{proof}

\section{Applications and Open Problems} \label{Applications}

One of the primary uses of Theorems \ref{ContinuousIsomorphism} and \ref{ContinuousZeroIsomorphism} is to study the continuous dual $\mathcal F_S'$ of $\mathcal F_S$.  Proposition \ref{FSubspace} provides the primary link 
between $\mathcal F_S'$ and $\lcc'$ or $\lcco'$.  
Let us first assume that $S$ is a proper subset of the set of all places of $\rat$.  In this case, Proposition \ref{FSubspace}\ref{SomeY} applies to show that $\mathcal F_S$
is a dense subset of $\lcc$ with respect to the norm \eqref{L1Norm}.  Every continuous linear functional on $\mathcal F_S$ extends to a continuous linear functional on $\lcc$, and therefore
\begin{equation} \label{FRewrite}
	\mathcal F_S' = \left\{ \Phi\in \lcc': \Phi(\mathcal F_S)\subseteq \rat\right\}.
\end{equation}
In view of Theorem \ref{ContinuousIsomorphism}, the following problem arises immediately.

\begin{question} \label{FQuestion}
	For which consistent maps $c\in \mathcal J_S^*$ do we have $\Phi_c(\mathcal F_S)\subseteq \rat$?
\end{question}

Note that we have phrased Question \ref{FQuestion} in terms of all consistent maps without reference to continuity.  However, if we wish to apply any answer to this question, it is likely that we need \eqref{FRewrite} in order
to do so.  This would require us to consider only those consistent maps that are associated to continuous linear functionals, i.e., those belonging to $\mathcal J_S'$.
Regardless, when $S$ is a set of non-Archimedean places of $\rat$, Question \ref{FQuestion} is completely resolvable.

\begin{thm} \label{nonArchRational}
	Let $S$ be a set of non-Archimedean places of $\rat$ and let $c\in \mathcal J_S^*$. Then $\Phi_c(\mathcal F_S)\subseteq \rat$ if and only if $c(K,v)\log p_v \in \rat$ for all $(K,v)\in \mathcal J_S$.
\end{thm}
\begin{proof}
	First assume that $c(K,v)\log p_v \in \rat$ for all $(K,v)\in \mathcal J_S$ and let $\alpha\in \alg$.  Suppose $K$ is a number field containing $\alpha$ so that 
	\begin{equation*}
		\Phi_c(f_{S,\alpha}) = \sum_{v\in M_{K,S}} c(K,v) \log \|\alpha\|_v.
	\end{equation*}
	We know that $\|\alpha\|_v$ is a rational power of $p_v$, so it follows that $\Phi_c(f_{S,\alpha})\in \rat$.
	
	Now assume that $\Phi_c(\mathcal F_S)\subseteq \rat$ and fix an element $(K,v)\in \mathcal J_S$.  By \cite[Lemma 3.1]{Samuels}, there exists $\alpha\in K$ such that $\|\alpha\|_v > 1$ but $\|\alpha\|_w = 1$ for all $w\in S\setminus \{v\}$.
	We may let $r$ be a positive rational number such that $\|\alpha\|_v = p_v^r$, so we obtain
	\begin{equation*}
		\Phi_c(f_{S,\alpha}) = c(K,v) \log \|\alpha\|_v = rc(K,v)\log p_v.
	\end{equation*}
	As we have assumed that $\Phi_c(f_{S,\alpha})\in \rat$, it follows that $c(K,v)\log p_v\in \rat$, as required.
\end{proof}

Let $\mathcal I_S^*$ be the $\rat$-vector space of consistent maps which satisfy the conditions of Theorem \ref{nonArchRational} and let $\mathcal I_S' = \mathcal I_S^* \cap \mathcal J_S'$.
Combined with Theorem \ref{ContinuousIsomorphism} and \eqref{FRewrite}, Theorem \ref{nonArchRational} yields an isomorphism from $\mathcal I_S'$ to $\mathcal F_S'$.  If $S$ is the set of all non-Archimedean places of $\rat$, then $\alpha\mapsto f_{S,\alpha}$
defines an isomorphism from $\mathcal V$ to $\mathcal F_S$.  Therefore, we reacquire \cite[Theorem 1.2]{Samuels} yielding an isomorphism from $\mathcal I_S'$ to $\mathcal V'$.  While Theorem \ref{ContinuousIsomorphism} is a 
more general result than \cite[Theorem 1.2]{Samuels}, it is not an entirely new proof.  Indeed, its proof extensively utilizes similar methods, and moreover, \cite[Lemma 3.1]{Samuels} is required in order to obtain \cite[Theorem 1.2]{Samuels}
as a corollary.

We close this section with a discussion of several open problems that relate to Question \ref{FQuestion} and Theorem \ref{nonArchRational}.
For the moment, it is unclear how one might apply Theorem \ref{nonArchRational} to obtain a new proof of the representation theorem \cite[Theorem 1.1]{Samuels} for $\mathcal V^*$.  As such, we find it a worthwhile problem
to investigate whether there is an analog of \eqref{FRewrite} for $\mathcal F_S^*$ and $\lcc^*$.   It is similarly unclear how to obtain an analog of Theorem \ref{nonArchRational}
when $S$ contains the Archimedean place of $\rat$.  Indeed, the proof of Theorem \ref{nonArchRational} requires one to determine the conditions under which $c(K,v)\log\|\alpha\|_v$ is rational for all $\alpha\in \alg$.
While this was fairly straightforward when $v$ was non-Archimedean, the situation is murkier when $v$ is Archimedean.

When considering the case where $S$ is the set of all places of $\rat$, Proposition \ref{FSubspace}\ref{AllY} shows that $\mathcal F_S$ is a dense subset of $\lcco$ with respect to the norm \eqref{L1Norm}.
As an analog of \eqref{FRewrite}, we quickly obtain that
\begin{equation} \label{FRewrite2}
	\mathcal F_S' = \left\{ \Phi\in \lcco': \Phi(\mathcal F_S)\subseteq \rat\right\},
\end{equation}
so Question \ref{FQuestion} is still valid in this case.  However, if we seek an analog of Theorem \ref{nonArchRational}, we face a similar challenge here as was described above.  Specifically, we are unable to describe the conditions
under which $c(K,v)\log\|\alpha\|_v$ is rational for all $\alpha\in \alg$.  We also note an additional challenge that arises from working in $\lcco$ rather than $\lcc$.  Our proof of Theorem \ref{nonArchRational} makes
extensive use of approximations to indicator functions in $\lcc$.  Unfortunately, those indicator functions do not belong to $\lcco$, and hence, some modification to the proof of Theorem \ref{nonArchRational} would be required.

\section{Further Generalizations of Consistent Maps} \label{Generalization}

As we noted at the end of the introduction, we could have taken $S$ to be a set of places of an arbitrary number field $F$ rather than restricting our definitions to $\rat$.  In fact, all of our proofs can be easily 
adapted to yield more general results.  However, there is little to be gained by generalizing our theorems in such a way.  In this section, we propose this alternate definition of consistency, and subsequently, show
it to be essentially equivalent to our existing definition.

Suppose that $F$ is a number field and $S$ is a nonempty set of places of $F$.  If $K$ is a finite extension of $F$, we write $M_{K,S}$ for the set of places of $K$ that divide a place in $S$.  Now define
\begin{equation*}
	\mathcal J(F,S) = \left\{(K,v):[K:F] < \infty,\ v\in M_{K,S}\right\}
\end{equation*}
and note the obvious fact that if $S\subseteq T$ then $\mathcal J(F,S)\subseteq \mathcal J(F,T)$.  In the special case $F=\rat$ considered throughout the majority of this article, we have $\mathcal J(\rat,S) = \mathcal J_S$.
 A map $c:\mathcal J(F,S)\to \real$ is called {\it consistent} if 
\begin{equation*}
	c(K,v) = \sum_{w\mid v} c(L,w)
\end{equation*}
for all $(K,v)\in \mathcal J(F,S)$ and all finite extensions $L/K$.  Just as before, we may define addition and scalar multiplication of consistent maps in the usual way.
Specifically, if $c,d:\mathcal J(F,S)\to \real$ are consistent maps and $r\in \real$ then we let
\begin{equation} %\label{ConsistentArithmetic}
	(c+d)(K,v) = c(K,v) + d(K,v)\quad\mbox{and}\quad (rc)(K,v) = rc(K,v).
\end{equation}
It is easily verified that these operations cause the set of consistent maps $c:\mathcal J(F,S)\to \real$ to be a vector space over $\real$, which we shall denote by $\mathcal J^*(F,S)$.
The following lemma shows that consistency remains intact when considering a finite extension of $F$.

\begin{lem} \label{ConsistentTower}
	Suppose that $F$ is a number field and $E$ is a finite extension of $F$.  If $S_E$ and $S_F$ are sets of places of $E$ and $F$, respectively, such that $S_E \subseteq M_{E,S_F}$, then
	$\mathcal J(E,S_E) \subseteq \mathcal J(F,S_F)$.  Moreover, if $c:\mathcal J(F,S_F)\to \real$ is a consistent map then the restriction of $c$ to $\mathcal J(E,S_E)$ is also consistent.
\end{lem}
\begin{proof}
	Let $(K,v)\in \mathcal J(E,S_E)$ so that $[K:E]<\infty$ and $v\in M_{K,S_E}$.  Since $E/F$ is assumed to be finite we certainly have that $[K:F] < \infty$.  Additionally, $v$ must divide a place
	$q\in S_E$.  We have assumed that $S_E \subseteq M_{E,S_F}$ so that $q$ must divide a place $p\in S_F$, and it follows that $v\mid p$.  We have now established that $(K,v)\in \mathcal J(F,S_F)$,
	and hence, $\mathcal J(E,S_E) \subseteq \mathcal J(F,S_F)$.
	
	Now suppose that $c:\mathcal J(F,S_F)\to \real$ is a consistent map and $d:\mathcal J(E,S_E)\to \real$ is its restriction to $\mathcal J(E,S_E)$.  To see that $d$ is consistent, let $(K,v)\in \mathcal J(E,S_E)$ and 
	let $L/K$ be a finite extension.  We have already shown that $(K,v)\in \mathcal J(F,S_F)$, so by consistency of $c$ we have that
	\begin{equation*}
		c(K,v) = \sum_{w\mid v} c(L,w).
	\end{equation*}
	Clearly $L/E$ is a finite extension and $w\in S_E$ so that $(L,w)\in \mathcal J(E,S_E)$.  It follows that $d(L,w)$ is well-defined and equals $c(L,w)$.  Hence
	\begin{equation*}
		d(K,v) = \sum_{w\mid v} d(L,w).
	\end{equation*}
	as required.
\end{proof}

We now find that our generalized definition of consistent is essentially independent of $F$.

\begin{thm} \label{ConsistentIsomorphism}
	Suppose that $S$ is a set of places of $\rat$, $F$ is a number field, and $T = M_{F,S}$.  Then the restriction map $\rho: \mathcal J^*(\rat,S)\to \mathcal J^*(F,T)$ given by
	$\rho(c) = c\mid_{\mathcal J(F,T)}$ is a vector space isomorphism.
\end{thm}
\begin{proof}
	It is straightforward to observe that $\rho$ is a linear map, so it remains only to show it to be a bijection.  We let $c:\mathcal J(F,T)\to \real$ be a consistent map and seek to show that
	$c$ has a unique extension to a consistent map $d:\mathcal J(\rat,S)\to \real$.
	
	Given a point $(K,v)\in \mathcal J(\rat,S)$, let $L$ be a number field containing both $K$ and $F$.  Following the notation of \cite{AllcockVaaler}, we write $W_v(L/K)$ for the set of places of $L$ that divide $v$.
	Also, let $p$ the unique place of $\rat$ for which $v\mid p$.  If $w\in W_v(L/K)$ then certainly $w\mid p$,
	and moreover, $w$ divides a place $q$ of $F$ such that $q\mid p$.  We conclude that $q\in M_{F,S} = T$ so that $(L,w)\in \mathcal J(F,T)$ for all $w\in W_v(L/K)$.  We have now shown that $c(L,w)$ is well
	defined for all $w\in W_v(L/K)$.
	
	Next, we set
	\begin{equation*}
		d_L(K,v) = \sum_{w\mid v} c(L,w)
	\end{equation*}
	and note that $d_L$ is independent of $L$.  Indeed, if $M$ is a finite extension of $L$ then the consistency of $c$ implies that
	\begin{equation*}
		d_L(K,v) = \sum_{w\mid v} \sum_{x\mid w} c(M,x) = \sum_{x\mid v} c(M,x) = d_M(K,v).
	\end{equation*}
	Therefore, we may define $d:\mathcal J(\rat,S)\to \real$ by
	\begin{equation} \label{ExtensionDefinition}
		d(K,v) = \sum_{w\mid v} c(L,w),
	\end{equation}
	where $L$ is any number field containing both $K$ and $F$.  If $(K,v)\in \mathcal J(F,T)$, then we may apply \eqref{ExtensionDefinition} with $L=K$ to see that $d$ is an extension of $c$.
	
	To see that $d$ is consistent, let $(K,v)\in \mathcal J(\rat,S)$ and $L/K$ a finite extension.  Further assume that $M$ is a finite extension of both $L$ and $F$.  Then using definition \eqref{ExtensionDefinition}, we obtain
	\begin{equation*}
		\sum_{w\mid v} d(L,w) = \sum_{w\mid v}\sum_{x\mid w} c(M,x) = \sum_{x\mid v} c(M,x) = d(K,v).
	\end{equation*}	
	
	Finally, we assume that $d_1:\mathcal J(\rat,S)\to \real$ and $d_2:\mathcal J(\rat,S)\to \real$ are consistent maps that extend $c$.  Again, let $(K,v)\in \mathcal J(\rat,S)$ and let $L$ be a finite extension of both $K$ and $F$.
	Then the consistency of $d_1$ and $d_2$ means that
	\begin{equation*}
		d_1(K,v) = \sum_{w\mid v} d_1(L,w).
	\end{equation*}
	However, $(L,w)\in \mathcal J(F,T)$, and since both $d_1$ and $d_2$ are extensions of $c$, we must have
	\begin{equation*}
		d_1(K,v) = \sum_{w\mid v} d_2(L,w) = d_2(K,v),
	\end{equation*}
	establishing that $d_1 = d_2$.  It now follows that $\rho$ is a bijection.
\end{proof}

We note that definition \eqref{ExtensionDefinition} actually provides a formula for $\rho^{-1}$.  Indeed, if $c\in \mathcal J(F,T)$ then $\rho^{-1}(c)$ satisfies the formula
\begin{equation*}
	[\rho^{-1}(c)](K,v) = \sum_{w\mid v} c(L,w),
\end{equation*}
where $L$ is any number field containing $F$ and $K$.


\begin{thebibliography}{10}

\bibitem{AkhtariVaaler}
S.~Akhtari and J.~D. Vaaler.
\newblock Heights, regulators and {S}chinzel's determinant inequality.
\newblock {\em Acta Arith.}, 172(3):285--298, 2016.

\bibitem{AllcockVaaler}
D.~Allcock and J.~D. Vaaler.
\newblock A {B}anach space determined by the {W}eil height.
\newblock {\em Acta Arith.}, 136(3):279--298, 2009.

\bibitem{Bogachev}
V.~I. Bogachev.
\newblock {\em Measure theory. {V}ol. {I}, {II}}.
\newblock Springer-Verlag, Berlin, 2007.

\bibitem{DunfordSchwartz}
N.~Dunford and J.~T. Schwartz.
\newblock {\em Linear operators. {P}art {I}}.
\newblock Wiley Classics Library. John Wiley \& Sons, Inc., New York, 1988.
\newblock General theory, With the assistance of William G. Bade and Robert G.
  Bartle, Reprint of the 1958 original, A Wiley-Interscience Publication.

\bibitem{FiliMinerNorms}
P.~Fili and Z.~Miner.
\newblock Norms extremal with respect to the {M}ahler measure.
\newblock {\em J. Number Theory}, 132(1):275--300, 2012.

\bibitem{FiliMinerOrthogonal}
P.~Fili and Z.~Miner.
\newblock Orthogonal decomposition of the space of algebraic numbers and
  {L}ehmer's problem.
\newblock {\em J. Number Theory}, 133(11):3941--3981, 2013.

\bibitem{FiliMinerDirichlet}
P.~Fili and Z.~Miner.
\newblock A generalization of {D}irichlet's unit theorem.
\newblock {\em Acta Arith.}, 162(4):355--368, 2014.

\bibitem{GrizzardVaaler}
R.~Grizzard and J.~D. Vaaler.
\newblock Multiplicative approximation by the {W}eil height.
\newblock {\em Trans. Amer. Math. Soc.}, 373(5):3235--3259, 2020.

\bibitem{Gubler}
W.~Gubler.
\newblock Heights of subvarieties over {$M$}-fields.
\newblock In {\em Arithmetic geometry ({C}ortona, 1994)}, volume XXXVII of {\em
  Sympos. Math.}, pages 190--227. Cambridge Univ. Press, Cambridge, 1997.

\bibitem{Hughes}
A.~M. Hughes.
\newblock A {W}eil {B}anach algebra for multiplicative algebraic numbers.
\newblock {\em J. Number Theory}, 196:244--271, 2019.

\bibitem{KellySamuels}
J.~P. Kelly and C.~L. Samuels.
\newblock Direct limits of ad\`ele rings and their completions.
\newblock {\em Rocky Mountain J. Math.}, 50(3):1021--1043, 2020.

\bibitem{Rudin}
W.~Rudin.
\newblock {\em Real and complex analysis}.
\newblock McGraw-Hill Book Co., New York, third edition, 1987.

\bibitem{Samuels}
C.~L. Samuels.
\newblock A classification of {$\Bbb Q$}-valued linear functionals on
  {$\overline {\Bbb Q}^{\times }$} modulo units.
\newblock {\em Acta Arith.}, 205(4):341--370, 2022.

\bibitem{Vaaler}
J.~D. Vaaler.
\newblock Heights on groups and small multiplicative dependencies.
\newblock {\em Trans. Amer. Math. Soc.}, 366(6):3295--3323, 2014.

\end{thebibliography}
\end{document}